\documentclass[10pt]{article}

\usepackage{graphicx}        
\usepackage{multicol}        
\usepackage[bottom]{footmisc}
\usepackage{amssymb,amsmath}

 \usepackage{authblk}

\makeindex             

\begin{document}

\title{3D-2D Stokes-Darcy coupling for the modelling of seepage with
 an application to fluid-structure interaction with contact}
\date{}
\author[$\star$]{Erik Burman}
\author[$\dagger$]{Miguel A. Fern\'andez}
\author[$\ddagger$]{Stefan Frei}
\author[$\dagger$]{Fannie M. Gerosa}
\affil[$\star$]{\footnotesize  Department of Mathematics, University College London, London, UK--WC1E  6BT, United Kingdom}
\affil[$\dagger$]{\footnotesize  Inria, 75012 Paris,  
Sorbonne Universit\'e \& CNRS, UMR 7598 LJLL, 75005 Paris.}
\affil[$\ddagger$]{\footnotesize Department of Mathematics \& Statistics, University of Konstanz, Germany.}
%
%
\maketitle

\begin{abstract}
In this note we introduce a mixed dimensional Stokes-Darcy coupling where a
  $d$ dimensional Stokes' flow is coupled to a Darcy model on the
  $d-1$ dimensional boundary of the domain. The porous layer introduces
  tangential creeping flow along the boundary and allows for the modelling
  of boundary flow due to surface roughness. This leads to a new
  model of flow in fracture networks with reservoirs in an
  impenetrable bulk matrix. Exploiting this modelling capability, we then formulate a
  fluid-structure interaction method with contact, where the porous layer allows
  for mechanically consistent contact and release. Physical seepage in the contact zone due
  to rough surfaces is modelled by the porous layer. Some numerical
  examples are reported, both on the Stokes'-Darcy coupling alone and on
  the fluid-structure interaction with contact in the porous boundary layer.
\end{abstract}

\section{Introduction}
In numerous environmental or biomedical applications there is a need to model the coupling between
a flow in a reservoir and flow in a surrounding porous
medium. This is particularly challenging if the porous medium is
fractured and the bulk matrix has very low permeability.
Typically the fractures are modelled as $d-1$ dimensional
manifolds, embedded in a $d$ dimensional porous bulk matrix. 
For the modelling of the fractured porous medium we refer to
\cite{martin-jaffre-robert-26}. Observe however that if the bulk
permeability is negligible the fluid in the reservoir can not penetrate
into the fractures since the $d-1$ dimensional manifolds have an
intersection of the reservoir boundary of $d-1$ measure zero. This
means that such a model can not be used for the fluid flow between two
reservoirs connected by a fracture in an impenetrable medium. Here we
propose to introduce a Darcy equation for the tangential flow on the
boundary of the reservoir. Since this equation is set on a $d-1$
dimensional manifold it can provide an interface allowing for flow
from the reservoir to the cracks. The flow on the boundary
communicates with the flow in the cracks through continuity of
pressure and conservation expressed by Kirchhoff's law. This gives a cheap and flexible
model for flow in reservoirs connected by fractures. 

Our original motivation for this model is the particular case of fluid
structure interaction with contact where the situation described above
occurs when two boundaries enter in contact provoking a change of
topology of the fluid domain. It has recently been observed by several
authors \cite{AgerWalletal,BurmanFernandezFrei} that the consistent modelling of fluid-structure
interaction with contact requires a fluid model, in particular a pressure, also in the contact
zone. Indeed, some seepage is expected to occur due to permeability of
the contacting bodies or their surface roughness. Otherwise there is no continuous mechanism for the release of
contact and non-physical voids can occur. For instance, it was argued
in \cite{AgerWalletal} that a consistent modelling of FSI with contact requires a
complete modelling of the FSI-poroelastic coupling. Similar ideas were
introduced in \cite{BurmanFernandezFrei}, but for computational
reasons. Indeed, in the latter reference an elastic body immersed in
a fluid enters in contact with a rigid wall and to allow for a
consistent numerical modelling the permeability of the wall is
relaxed. This motivates the introduction of an artificial porous medium
whose permeability goes to zero with the mesh-size. Both approaches allow for the
seepage that appears to be necessary for physical contact and
release. 
However, in case the contacting solids are (modelled as) impenetrable, this seepage
must be due to porous media flow in a thin layer in the contact zone due to surface roughness.
The complete modelling of the poroelastic interaction of \cite{AgerWalletal} or the
bulk porous medium flow of \cite{BurmanFernandezFrei} then appears
artificial and unnecessarily expensive. For such situations the mixed
dimensional modelling suggested above can offer an attractive
compromise between model detail and computational cost.

In this note, we will focus exclusively on the modelling aspect. The
coupled Stokes-Darcy model is introduced in section
\ref{sec:stda}. Then, in section \ref{sec.FSIcontact}, we show how the ideas of
\cite{BurmanFernandezFrei} can be used to model FSI with contact
together with the mixed-dimensional fluid system. Finally, we
illustrate the two model situations numerically in section
\ref{sec:num}. First, the Stokes'-Darcy reservoir coupling (section
\ref{sec:numstda})  and then the
full FSI with contact (section \ref{sec:numFSIcon}). In the latter
case, we also give comparisons with the results from
\cite{BurmanFernandezFrei}. 
The numerical analysis of the resulting methods will be the subject of
future work.
\section{The coupled Stokes-Darcy system} \label{sec:stda}

We consider the coupling of a Darcy system in a thin-walled domain $\Omega_l = \Sigma_l \times (-\frac{\epsilon}{2},\frac{\epsilon}{2}) \in \mathbb{R}^d$ for $d=2,3$
with a Stokes equation in the bulk domain $\Omega_f$. The Darcy problem on $\Omega_l$ writes
\begin{equation}\label{eq:darcy}
\left\{
\begin{aligned}
 u_l+ K \nabla p_l = 0 & \\
\nabla \cdot u_l = 0 & \\
\end{aligned}
\right.\quad\mbox{in}\quad \Omega_l,
\end{equation}
where $u_l$ denotes the Darcy velocity, $p_l$ the Darcy pressure and $K$ is a $d\times d$ matrix that allows for the decomposition
\begin{align*}
K\nabla p_l = K_\tau \nabla_{\tau} p_l + K_n \partial_n p_l.
\end{align*}
We denote the upper boundary of $\Omega_l$ which couples to $\Omega_f$ by $\gamma_f$ and 
the outer boundary by $\gamma_o$. The normal vector $n$ of the middle surface $\Sigma_l$ of $\Omega_l$ is chosen in such a way that it points towards $\gamma_o$.

By averaging across the thickness $\epsilon$, Martin, Jaffr\'e and Roberts derived in \cite{martin-jaffre-robert-26} an effective equation 
for the averaged pressure across the thickness
 $$
 P_l := \frac{1}{\epsilon}\int_{-\frac{\epsilon}{2}}^{\frac{\epsilon}{2} }p_l .
 $$
 Under the modelling assumption that the average pressure is equal to the mean of the pressures on the upper and lower boundary
 \begin{equation}\label{eq:meanp}
P_l= \frac{1}{2}\left(p_l\vert{\gamma_f} + p_l\vert{\gamma_o}  \right)  \quad\mbox{in}\quad \Sigma_l,
\end{equation}
the authors derived the system 
\begin{equation}\label{eq:PP}
\left\{
\begin{aligned}
- \nabla_\tau \cdot \left( \epsilon K_\tau \nabla_\tau P_l  \right) &=  u_{l,n}\vert{\gamma_f}  - u_{l,n}\vert_{\gamma_o} \\
 p_l\vert{\gamma_f} &= P_l  + \frac{\epsilon K_n^{-1} }{4}\left( u_{l,n}\vert{\gamma_o} + u_{l,n}\vert{\gamma_f}   \right)  
 \end{aligned}\right. \quad\mbox{in}\quad \Sigma_l.
\end{equation}
Here, $u_{l,n}=u_l\cdot n$ denotes the normal component of the velocity and $\tau$ is a tangential vector of $\Sigma_l$.
We will couple \eqref{eq:PP} to Stokes flow in $\Omega_f$ 
\begin{equation}
\left\{
\begin{aligned}
\rho_f \partial_t u_f- \nabla \cdot \sigma_f (u_f,  p_f) = 0 &\\
\nabla \cdot u_f = 0 &
\end{aligned}
\right. \quad\mbox{in}\quad \Omega_f,
\end{equation}
where $u_f$ denotes the fluid velocity,  $p_f$ the pressure, $\rho^{\rm f}$ the fluid density, 
$$
\sigma_f (u_f,p_f) := \mu ( \nabla u_f + \nabla u_f^T) - p_f I,
$$
the fluid Cauchy stress tensor and $\mu$ the dynamic viscosity. 
  We
assume that the coupling to the Darcy system \eqref{eq:darcy} on $\gamma_f$ takes place 
 via the interface conditions 
 \begin{equation}\label{eq:intf}
 \left\{
 \begin{aligned}
 \sigma_{f,nn} = - p_l  &  \\
 \tau^T \sigma_{f} n = 0 &  \\
 u_{f,n} = u_{l,n} & \\
 \end{aligned}\right. \quad\mbox{on}\quad \gamma_f,
 \end{equation}
 where $\sigma_f=\nabla u_f -p_f I$ and $\sigma_{f,nn}=n^T\sigma_f n$.
In the lower porous wall $\gamma_o$ we assume for simplicity that $u_{l,n} = 0$.
Then, the relations \eqref{eq:PP} can be written as 
$$
\left\{
\begin{aligned}
- \nabla_\tau \cdot \left( \epsilon K_\tau \nabla_\tau P_l  \right) =    u_{f,n}  & \\
 \sigma_{f,nn} = - P_l  - \frac{\epsilon K_n^{-1} }{4} u_{f,n}   & 
\end{aligned}
\right. \quad\mbox{in}\quad \Sigma_l.
$$
Note that the only remaining porous medium variable is the averaged pressure $P_l$.
In the limit of permeability $K_n \to 0$, the system converges to a pure Stokes system with slip conditions on $\gamma_f$
with an extension of the fluid forces into the porous medium pressure $P_l$.

We have the following coupled variational problem for $(u_f,p_f,P_l)$: 
\begin{equation}
\left\{
\begin{aligned}
\rho_f  (\partial_t u_f,v_f)_{\Omega_f}  +   (\sigma_f (u_f,p_f), \nabla v_f)_{\Omega_f} + (q_f,\nabla \cdot u_f)_{\Omega_f}& \\
 + \big(P_l,v_{f,n} \big)_{\Sigma_l} + \frac{\epsilon K_n^{-1}}{4}\big( u_{f,n},v_ {f,n}   \big)_{\Sigma_l} &= 0,\\
(\epsilon K_\tau \nabla_\tau P_l, \nabla_{\tau} q_l)_{\Sigma_l} - ( u_{f,n}, q_l \big)_{\Sigma_l} &=0,
\end{aligned}
\right.
\end{equation}
for all $v_f,q_f,q_l$, where $n=n_f$ is the outer normal of the fluid domain $\Omega_f$.

\section{The fluid-structure-poroelastic-contact interaction system} 
\label{sec.FSIcontact}

Now, we consider a fluid-structure-contact interaction system with a thin porous layer on the part of the exterior boundary, 
where contact might take place. 
The moving boundary of the solid is denoted by $\Sigma(t)$ and the porous layer by $\Sigma_l$. In absence of contact, we have 
the following system of equations
\begin{equation*}
\begin{aligned}
& \left\{
\begin{aligned}
\rho_f \partial_t u_f  -\nabla \cdot \sigma_f(u_f,p_f)= 0 &\\
\nabla \cdot u_f = 0 &
\end{aligned}
\right.\quad\mbox{in}\quad \Omega_f(t),\\
& \rho_s \partial_t \dot{d} -\nabla \cdot \sigma_s(d) = 0  \quad\mbox{in}\quad \Omega_s(t),\\
& u_f = \dot{d}, \quad \sigma_s n =\sigma_f n  \quad\mbox{in}\quad \Sigma(t),
\end{aligned}
\end{equation*}
\begin{equation}
\left\{
\begin{aligned}
- \nabla_\tau \cdot \left( \epsilon K_\tau \nabla_\tau P_l  \right) =   u_{l,n}|_{\gamma_f} \label{ulnepsK}\\   
 \sigma_{f,nn} = \underbrace{- P_l  - \frac{\epsilon K_n^{-1} }{4} u_{l,n}|_{\gamma_f}}_{\sigma_p} \\
 \tau^T \sigma_f n =0 &
\end{aligned}
\right.  \quad\mbox{in}\quad \Sigma_l,
\end{equation}
where, in addition to the quantities introduced above, $\rho_s$ denotes the solid density, 
$d$ stands for the solid displacement and $\sigma_s$ denotes the tensor of linear elasticity
\begin{align*}
 \sigma_s = \frac{\lambda_s}{2} \left(\nabla d + \nabla d^T\right) + \frac{\mu_s}{2}\text{tr}\left(\nabla d + \nabla d^T\right).
\end{align*}
In addition, we impose that the solid $\Omega_s$ can not penetrate into the porous medium $\Sigma_l$
\begin{equation}\label{eq:ContCond}
 d_n - g \leq 0, \quad \lambda \leq 0, \quad \lambda (d_n - g)  = 0 \quad \text{ on } \Sigma(t).
\end{equation}
Here, $g$ denotes the gap function to $\Sigma_l$ and $\lambda$ is a Lagrange multiplier for the no-penetration condition defined by
\begin{align*}
 \lambda &= \sigma_{s,nn} - \sigma_{f,nn} \qquad\text{ on } \Sigma(t)\setminus \Sigma_l, \\
 \lambda &= \sigma_{s,nn} - \sigma_p \qquad\quad \text{ on } \Sigma(t)\cap \Sigma_l.
\end{align*}
The ``switch'' on the right-hand side occurs, as the solid on one side of $\Sigma(t)$ couples either to the fluid $\Omega_f$ or the porous
medium $\Sigma_l$ on the other side of $\Sigma(t)$. The conditions~\eqref{eq:ContCond} can equivalently be written as
\begin{align*}
 \lambda = \gamma_C \big[  \underbrace{d_n - g -\gamma_C^{-1} \lambda}_{P_\gamma} \big]_+ \quad \text{ on } \Sigma(t)
\end{align*}
for arbitrary $\gamma_C>0$. Using this notation, we can characterise the zone of ``active'' contact as follows
$$
\Sigma_c(t) = \left\{ x \in \Sigma_(t) \, | \,  P_{\gamma} > 0  \right\}. 
$$

To summarise, we have the following interface conditions:
\begin{itemize}
\item Contact condition on $\Sigma(t)$:
$$
d_n - g \leq 0, \quad \lambda \leq 0, \quad \lambda (d_n - g)  = 0 \quad \text{on} \quad \Sigma(t).
$$
\item Kinematic coupling on $\Sigma_{fsi}(t) = \Sigma(t) \backslash \Sigma_l$ 
$$
u_f = \dot{d}\quad \mbox{on}\quad \Sigma_{fsi}(t).
$$
\item Dynamic coupling on $\Sigma(t)$:
\begin{align*}
\sigma_sn  &= \lambda n - \sigma_{p} n = \gamma_C [P_\gamma]_+ n - \sigma_{p}n \quad \text{on }\Sigma(t) \cap \Sigma_l,\\
\sigma_sn  &= \lambda n - \sigma_f n = \gamma_C [P_\gamma]_+ n - \sigma_f n \quad \text{on }\Sigma(t) \setminus \Sigma_l.
\end{align*}
\end{itemize}
\vspace{0.3cm}

\noindent We have the following Nitsche-based variational formulation: \textit{Find } $u_f\in {\mathcal V}_f, p_f\in {\mathcal L}_f, d\in {\mathcal V}_s, P_l \in {\mathcal V}_l$ such that 
\begin{multline*}
(\partial_t u_f, v)_{\Omega_f} + (\partial_t \dot{d}, w)_{\Omega_s} + a_f\big(u_f,p_f; v,q\big) + a_s(d,w) \\
- (\sigma_fn,  v-w)_{\Sigma(t)\setminus \Sigma_l} 
- (u_f-\dot{d}, \sigma_f(v,-q))_{\Sigma(t)\setminus \Sigma_l}
+ \frac{\gamma_{\text{fsi}}}{h}(u_f-\dot{d}, v-w)_{\Sigma(t)\setminus \Sigma_l} \\
- (\sigma_p,  v\cdot n)_{\Sigma_l\setminus \Sigma(t)}  
- (\sigma_p,  w\cdot n)_{\Sigma_l \cap \Sigma(t)}
+ \big([P_\gamma]_+, w\cdot n\big)_{\Sigma(t)}\\
+ (\epsilon K_\tau \nabla_\tau P_l, \nabla_\tau q_l)_{\Sigma_l}
- \big( u_{f,n}, q_l\big)_{\Sigma_l\setminus \Sigma(t)}  - \big(\dot{d}_n , q_l \big)_{\Sigma_l\cap\Sigma(t)} 
=0 \\
\forall v\in {\mathcal V}_f, q\in {\mathcal L}_f, w\in {\mathcal V}_s, q_l \in {\mathcal V}_l.
\end{multline*}
The porous stress $\sigma_p$ is given by
\begin{align}\label{sigmaP}
 \sigma_p=- P_l  + \frac{\epsilon K_n^{-1} }{4} u_{l,n}|_{\gamma_f} = \begin{cases} 
                                                              &- P_l  + \frac{\epsilon K_n^{-1} }{4} u_{f,n} \quad \text{ on } \Sigma_l \setminus \Sigma(t)\\
                                                              &- P_l  + \frac{\epsilon K_n^{-1} }{4} \dot{d}_n \quad \text{ on }\Sigma_l \cap \Sigma(t).
                                                             \end{cases}
\end{align}



\section{Numerical experiments}\label{sec:num}
Here we will report on some numerical experiments using the above
models. First we consider the mixed dimensional Stokes'-Darcy system
and then the fluid-structure interaction system with contact and porous layer
in the contact zone.

\subsection{Stokes-Darcy example}\label{sec:numstda}

\begin{figure}[h!]
\centering
\includegraphics[width=0.45\textwidth]{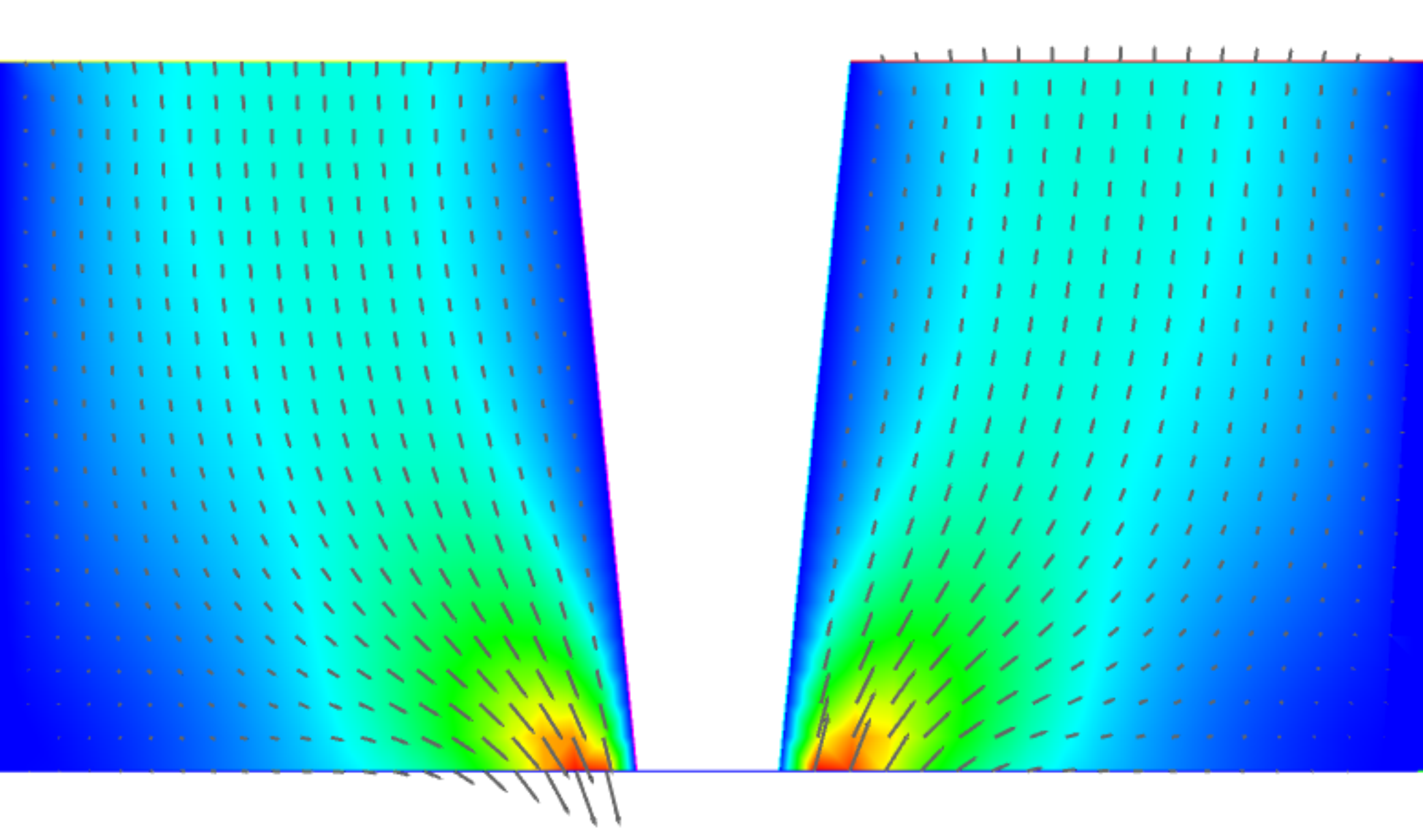}\hspace{1cm}
\includegraphics[width=0.45\textwidth]{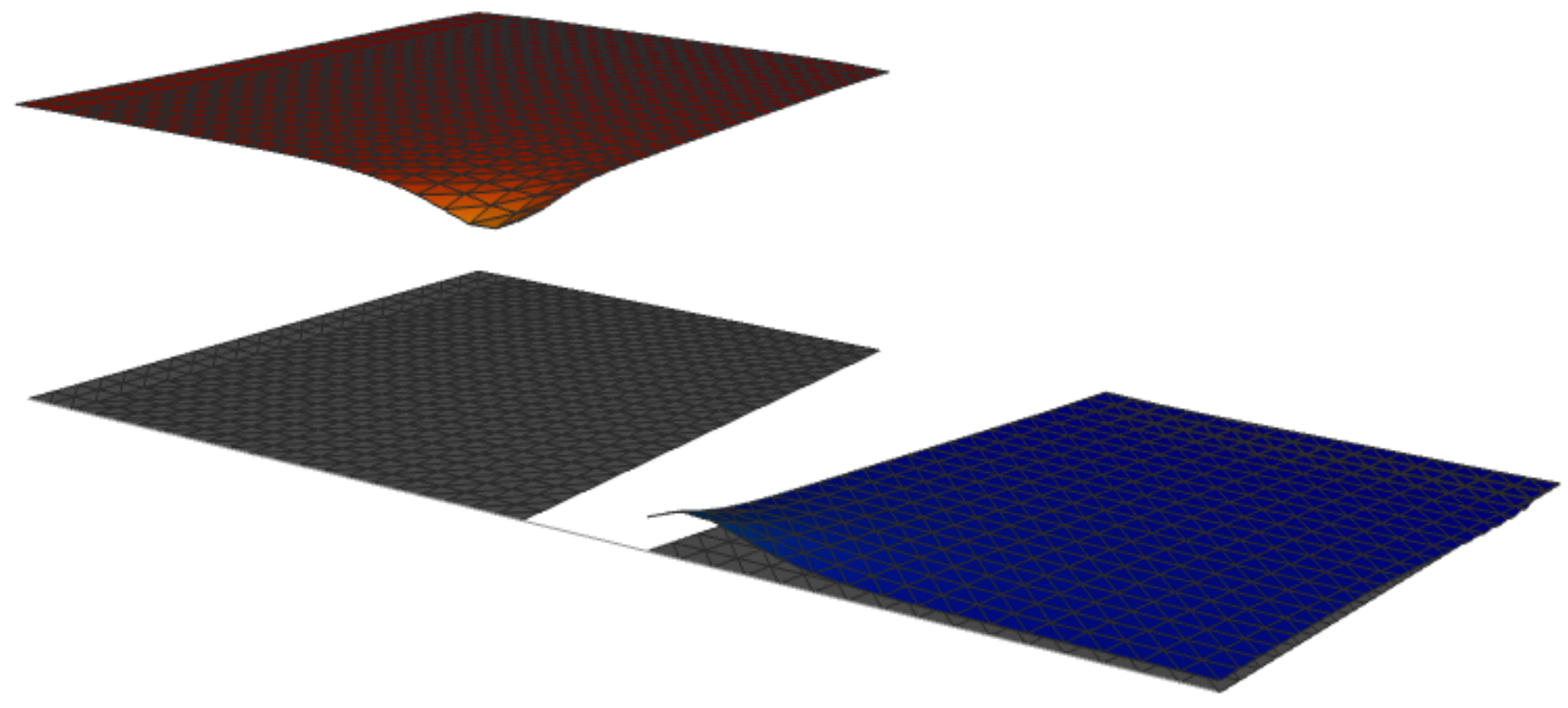}\\[10mm]
\includegraphics[width=0.5\textwidth]{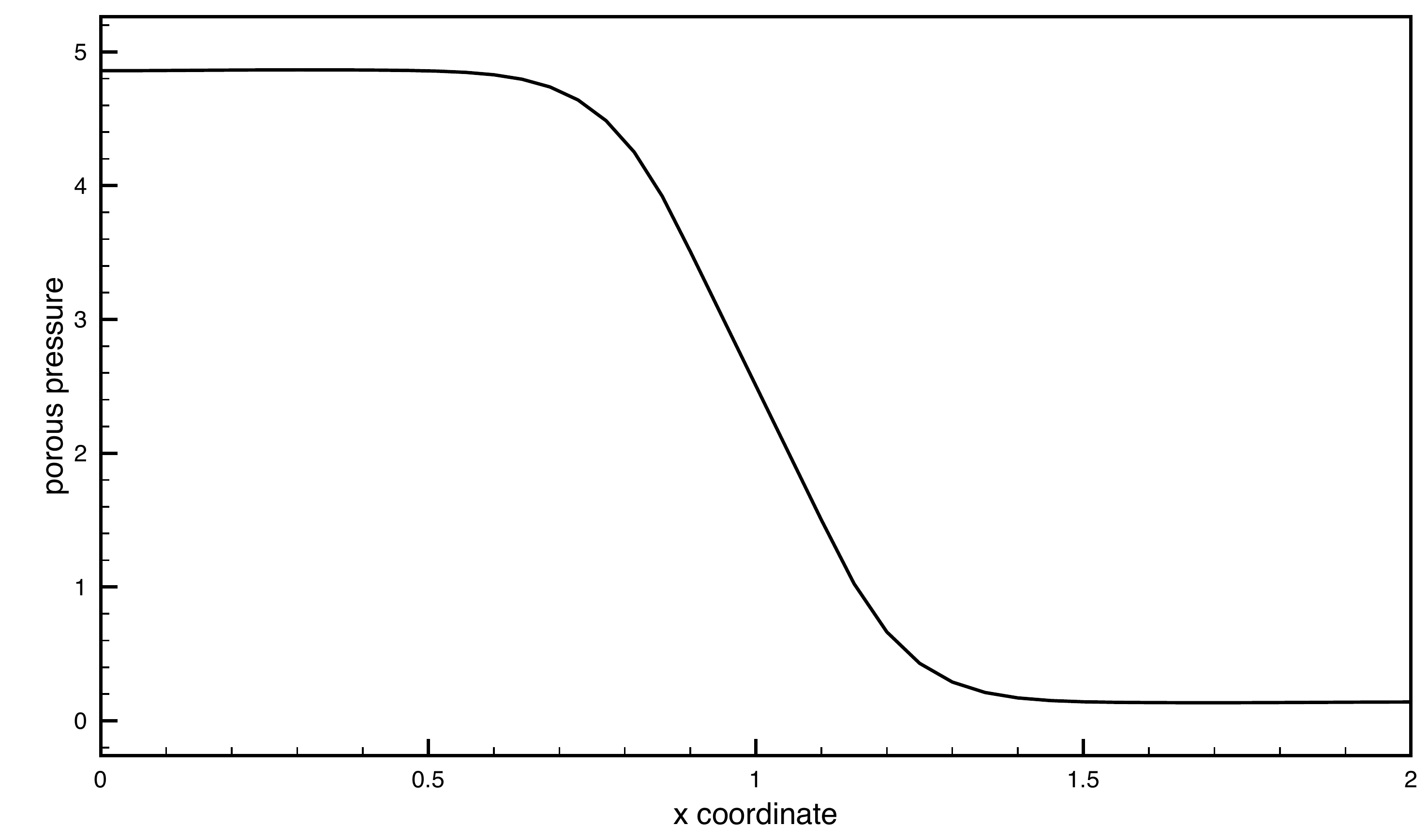}
\caption{\label{fig.nofsi-sim} 
Top left: Snapshot of the fluid velocity. 
Top right: Elevation of the fluid pressure. 
Bottom: Porous pressure.}
\end{figure}

In this example, we consider two disconnected fluid reservoirs, the domain $\Omega_f$, 
\begin{figure}[h!]
\centering
\includegraphics[width=0.8\textwidth]{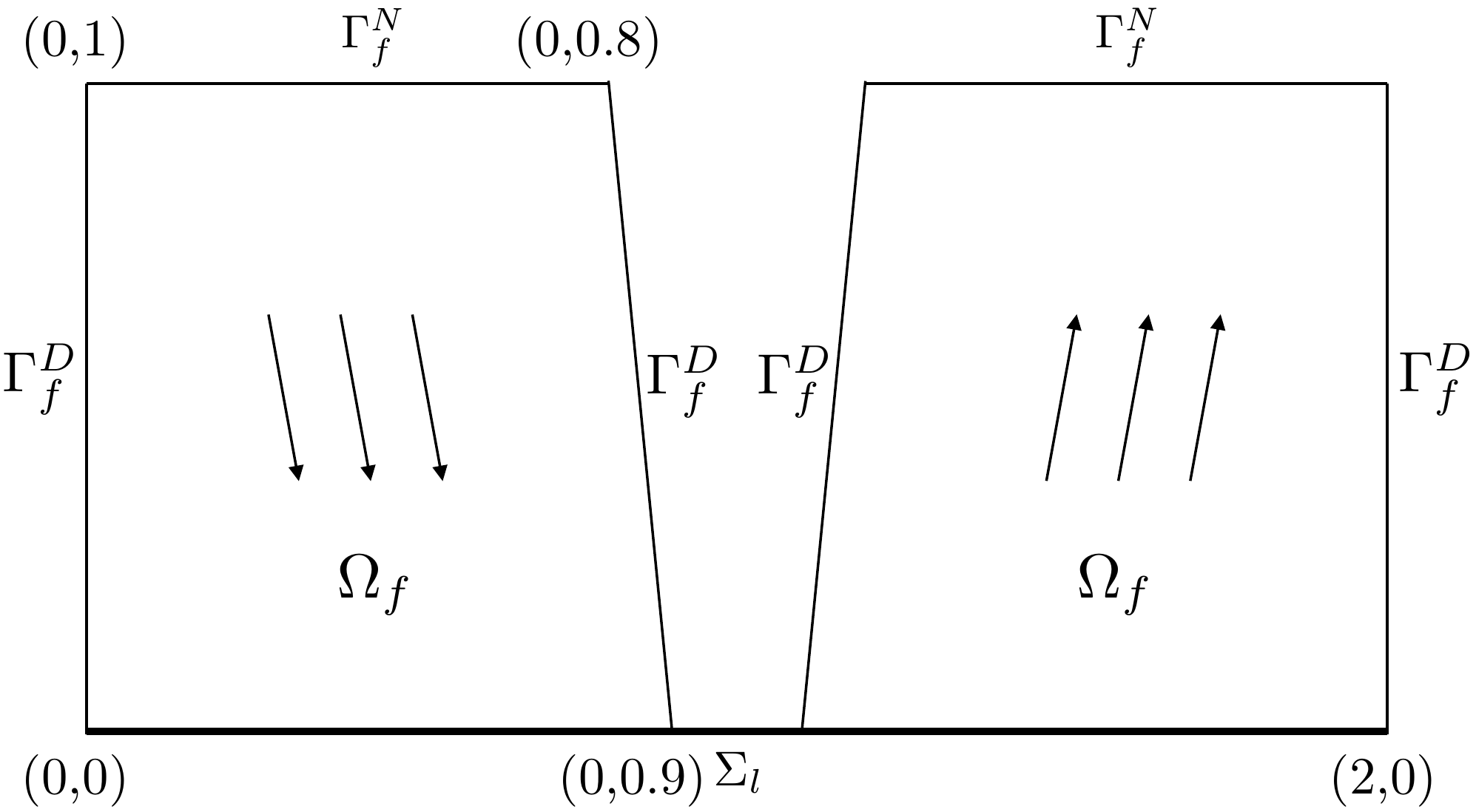}
\caption{\label{fig.nofsi} Geometrical configuration for the Stokes model with a thin-walled porous medium on the bottom wall.}
\end{figure}
connected through a thin-walled porous media located on the bottom wall $\Sigma_l$, as shown in 
Figure~\ref{fig.nofsi}.  The physical parameters are $\mu=0.03$, $\rho_f=1$, $\epsilon = 0.01$ and $K_\tau=K_n = 1$. 
We impose a pressure drop across the two parts of the boundary $\Gamma_f^N$. The purpose of this example is to illustrate how the 
porous model is able to connect the fluid flow between the two containers. This can be clearly inferred from the results reported in 
Figure~\ref{fig.nofsi-sim}, which respectively show a snapshot of the fluid velocity, the elevation of the fluid pressure and the associated porous pressure.

\subsection{Fluid-structure interaction with contact}\label{sec:numFSIcon}

\begin{figure}
\centering
\includegraphics[width=0.75\textwidth]{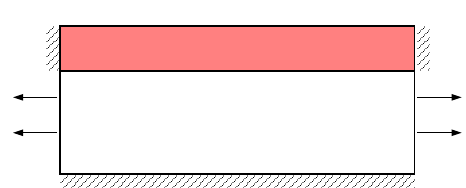}
\caption{\label{fig.FSIconf}Geometrical configuration for the FSI-contact model. We apply a porous medium model on the (rigid) lower wall, where contact might take place.}
\end{figure}

To test the FSI-contact model, we consider flow in a 2-dimensional pipe, where the upper wall is elastic, see Figure~\ref{fig.FSIconf}. Due to the application of a 
large pressure $\overline{P}$ on the left and right boundary, the
upper wall is deflected downwards until it reaches the bottom. Note
that when contact occurs, the configuration is topologically
equivalent to the situation in section \ref{sec:numstda}.
 Shortly before the time of impact we set $\overline{P}$ to zero, such that contact is realeased again after a certain time. This model 
 problem is taken from \cite{BurmanFernandezFrei}, where further details on the configuration and the discretisation can be found. To deal with the topology change in the fluid domain at the impact, 
we apply a \textit{Fully Eulerian} approach for the FSI problem~\cite{FreiPhD}.
In order to obtain a continuous and physically relevant transition from FSI to solid-solid contact, we use the FSI-contact model derived in section~\ref{sec.FSIcontact} 
and place a thin porous domain $\Sigma_l$ on the lower boundary.

\begin{figure}
\includegraphics[width=0.5\textwidth]{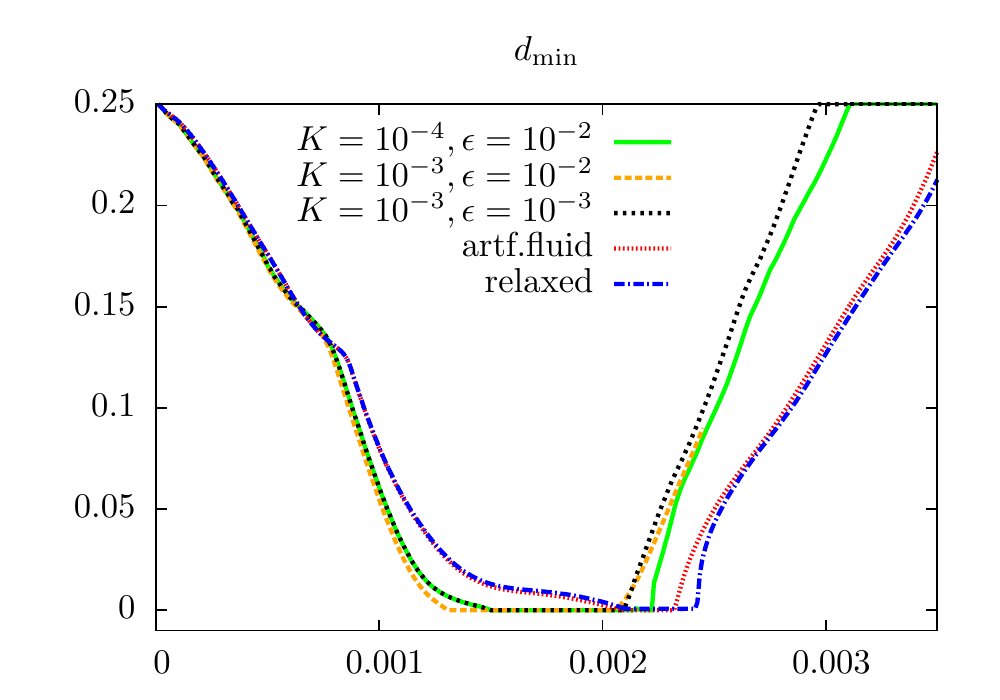}
\includegraphics[width=0.5\textwidth]{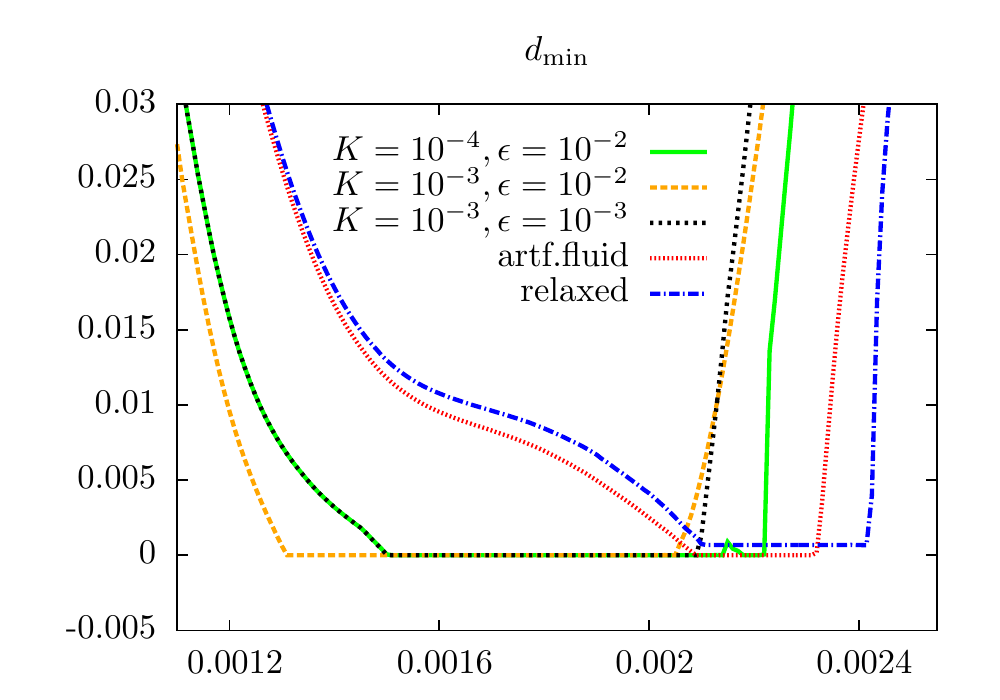}
  \caption{\label{fig.dmin} Minimal distance of $\Omega_s$ to the lower wall $\Sigma_p$ over time. Right: zoom-in around the contact interval. We 
 compare the new approach presented in Section~\ref{sec.FSIcontact} for different parameters with the \textit{artificial fluid} and the \textit{relaxed} contact 
  approach studied in~\cite{BurmanFernandezFrei}.
  }
\end{figure}

In Figure~\ref{fig.dmin} we compare this model for different parameters $K=K_\tau=K_n$ and $\epsilon$ with the approaches for FSI-contact problems 
introduced in \cite{BurmanFernandezFrei} in terms of the 
minimal distance of the solid to $\Sigma_p$ over time. In \cite{BurmanFernandezFrei} two approaches
were presented in order to extend the fluid stresses to the contact region during solid-solid contact, namely a so-called \textit{relaxed} and an 
\textit{artificial fluid} approach. It was observed that for the \textit{artificial fluid} approach contact happens earlier, as penetration of the 
fluid flow into the artificial region is prevented only asymptotically, i.e. $u_{f,n}\to 0\, (h\to 0)$ on $\Sigma_p$, in contrast to $u_{f,n}=0$ 
for the relaxed approach. In the model presented here, we have
similarly from \eqref{ulnepsK} and $u_{l,n}=u_{f,n}$ on $\Sigma_p$
\begin{align*}
 u_{f,n}=-\nabla_\tau\cdot (\epsilon K_\tau \partial_\tau P_l) \to 0 \quad (\epsilon K_{\tau} \to 0).
\end{align*}
For this reason we observe in Figure~\ref{fig.dmin} that the impact happens earlier for a larger value of $\epsilon K_\tau$.
The time of the release seems to depend also
on  $\epsilon K_n^{-1}$, which appears in the definition of $\sigma_p$ \eqref{sigmaP}. A detailed investigation of this dependence and the 
investigation of stability and convergence of the numerical method are subject to future work.

\section*{Acknowledgments}
Erik Burman was partially supported by the EPSRC grant: EP/P01576X/1. Stefan Frei acknowledges support by the DFG Research Scholarship FR3935/1-1.

 \bibliographystyle{plain}
 \bibliography{lit}

\begin{thebibliography}{1}

\bibitem{AgerWalletal}
C~Ager, B~Schott, AT~Vuong, A~Popp, and WA~Wall.
\newblock A consistent approach for fluid-structure-contact interaction based
  on a porous flow model for rough surface contact.
\newblock {\em Int J Numer Methods Eng}, 119(13):1345--1378, 2019.

\bibitem{FreiPhD}
S~Frei.
\newblock {\em {E}ulerian finite element methods for interface problems and
  fluid-structure interactions}.
\newblock PhD thesis, Heidelberg University, 2016.
\newblock http://www.ub.uni-heidelberg.de/archiv/21590.

\bibitem{martin-jaffre-robert-26}
V~Martin, J~Jaffr\'{e}, and JE~Roberts.
\newblock Modeling fractures and barriers as interfaces for flow in porous
  media.
\newblock {\em SIAM J. Sci. Comput.}, 26(5):1667--1691, 2005.

\bibitem{BurmanFernandezFrei}
E~Burman S~Frei and MA~Fern\'andez.
\newblock Nitsche-based formulation for fluid-structure interactions with
  contact.
\newblock {\em ESAIM: M2AN (published online)}.
\newblock https://doi.org/10.1051/m2an/2019072.

\end{thebibliography}

\end{document}